\documentclass[twoside,leqno,10pt, A4]{amsart}
\usepackage{amsfonts}
\usepackage{amsmath}
\usepackage{amscd}
\usepackage{amssymb}
\usepackage{amsthm}
\usepackage{amsrefs}
\usepackage{latexsym}
\usepackage{mathrsfs}
\usepackage{bbm}
\usepackage{enumerate}
\usepackage{color}
\begin{document}

\newtheorem{theorem}[subsection]{Theorem}
\newtheorem{proposition}[subsection]{Proposition}
\newtheorem{lemma}[subsection]{Lemma}
\newtheorem{corollary}[subsection]{Corollary}
\newtheorem{conjecture}[subsection]{Conjecture}
\newtheorem{prop}[subsection]{Proposition}
\numberwithin{equation}{section}
\newcommand{\mr}{\ensuremath{\mathbb R}}
\newcommand{\mc}{\ensuremath{\mathbb C}}
\newcommand{\dif}{\mathrm{d}}
\newcommand{\intz}{\mathbb{Z}}
\newcommand{\ratq}{\mathbb{Q}}
\newcommand{\natn}{\mathbb{N}}
\newcommand{\comc}{\mathbb{C}}
\newcommand{\rear}{\mathbb{R}}
\newcommand{\prip}{\mathbb{P}}
\newcommand{\uph}{\mathbb{H}}
\newcommand{\fief}{\mathbb{F}}
\newcommand{\majorarc}{\mathfrak{M}}
\newcommand{\minorarc}{\mathfrak{m}}
\newcommand{\sings}{\mathfrak{S}}
\newcommand{\fA}{\ensuremath{\mathfrak A}}
\newcommand{\mn}{\ensuremath{\mathbb N}}
\newcommand{\mq}{\ensuremath{\mathbb Q}}
\newcommand{\half}{\tfrac{1}{2}}
\newcommand{\f}{f\times \chi}
\newcommand{\summ}{\mathop{{\sum}^{\star}}}
\newcommand{\chiq}{\chi \bmod q}
\newcommand{\chidb}{\chi \bmod db}
\newcommand{\chid}{\chi \bmod d}
\newcommand{\sym}{\text{sym}^2}
\newcommand{\hhalf}{\tfrac{1}{2}}
\newcommand{\sumstar}{\sideset{}{^*}\sum}
\newcommand{\sumprime}{\sideset{}{'}\sum}
\newcommand{\sumprimeprime}{\sideset{}{''}\sum}
\newcommand{\shortmod}{\ensuremath{\negthickspace \negthickspace \negthickspace \pmod}}
\newcommand{\V}{V\left(\frac{nm}{q^2}\right)}
\newcommand{\sumi}{\mathop{{\sum}^{\dagger}}}
\newcommand{\mz}{\ensuremath{\mathbb Z}}
\newcommand{\leg}[2]{\Big(\frac{#1}{#2}\Big)}
\newcommand{\muK}{\mu_{\omega}}
\newcommand{\seteq}{:=}
\newcommand{\odd}{\mathrm{\ primary}}

\title[Murmurations of quadratic Hecke $L$-functions of the Gaussian field]{Murmurations of quadratic Hecke $L$-functions of the Gaussian field}

\author[P. Gao]{Peng Gao}
\address{School of Mathematical Sciences, Beihang University, Beijing 100191, China}
\email{penggao@buaa.edu.cn}

\author[L. Zhao]{Liangyi Zhao}
\address{School of Mathematics and Statistics, University of New South Wales, Sydney NSW 2052, Australia}
\email{l.zhao@unsw.edu.au}

\begin{abstract}
In this paper, we evaluate the murmuration density for the family of quadratic Hecke $L$-functions of the Gaussian field under the generalized Riemann hypothesis. 
\end{abstract}

\maketitle



\section{Introduction}

In \cite{HLOP25}, Y.-H. He, K.-H. Lee, T. Oliver and A. Pozdnyakov discovered a remarkable oscillation pattern in
the averages of Frobenius traces for certain families of elliptic curves, based on machine learning algorithms. This phenomenon is termed murmuration, because of its visual similarity to swarming patterns of flocks of birds.  N. Zubrilina \cite{Zubrilina25} provided the first theoretical confirmation (note that the preprint of \cite{Zubrilina25} appeared on the arxiv in 2023) on the existence of murmurations in the family of self-dual holomorphic newforms.  The results in \cite{Zubrilina25} consider a short average over conductors and illustrate correlations for an individual prime.  In \cite{Sarnak23}, P. Sarnak introduced a notion of murmuration density which takes additional averaging over primes into account.  It is also pointed out in \cite{Sarnak23} that there is a close relationship between murmurations and the one-level density of low-lying zeros of families of $L$-functions (see \cites{KS1, K&S}), which involves evaluating sums over a family of $L$-functions over analytic conductors up to $N$ and over the primes up to the scale of $N^{\theta}$ for some real $\theta>0$. There is a phase transition at some $\theta=\theta_0$ and the murmuration density interpolates the phase transition range.  In \cite{LOP25}, K.-H. Lee, T. Oliver and A. Pozdnyakov calculated murmuration densities for two families of Dirichlet characters, including one of quadratic Dirichlet characters under the generalized Riemann hypothesis (GRH). It was also shown that for this family, the corresponding murmuration density reveals the phase transition range from the one-level density case. \newline
  
The authors \cite{G&Zhao4} studied the one-level density of low-lying zeros of the family of quadratic Hecke $L$-functions in the Gaussian field.  The aim of this paper to apply the techniques developed there to evaluate the murmuration density for the corresponding family under GRH.  Let $K=\mq(i)$ be the Gaussian field and $\mathcal O_K=\mz[i]$ the ring of integers in $K$ throughout the paper.  We denote by $\zeta_K(s)$ the Dedekind zeta function of $K$ and $U_K=\{\pm 1, \pm i \}$ the group of units of $\mathcal O_K$ and $N(c)$ the norm of any $c \in K$. An element $c \in \mathcal O_K$ is said to be primary if $c \equiv 1 \pmod {(1+i)^3}$ and we recall that (see \cite[Chap. 9, Lemma 7]{I&R}) every ideal co-prime to $2$ in $\mathcal O_K$ has a unique primary generator. \newline

  We reserve the letter $\varpi$ for a prime in $\mathcal O_K$, i.e. the ideal generated by $\varpi$ is a prime ideal. We say an element $c \in \mathcal O_K$ is odd if $(c,2)=1$ and we say $c$ is square-free if $c$ is not divisible by the square of any prime.  Let $\mu_{[i]}$ denote the M\"obius function on $\mathcal O_K$ so that $\mu_{[i]}(c)=1$ for $c \in U_K$, $\mu_{[i]}(c)=(-1)^j$ if the ideal $(c)$ is the product of $j$ distinct prime ideals and $\mu_{[i]}(c)=0$ is $c$ is not square-free. \newline

Suppose $\Phi(t)$ is a non-negative smooth function with compact support on the set of positive real numbers. Define for $y$, $\delta>0$,
\begin{align*}
\begin{split}
M_{\Phi}(y,X,\delta)=\frac{\log X}{X^{1+\delta}}\sum_{\substack{N(\varpi) \in[yX,yX+X^{\delta}] \\ \varpi \odd }}\sum_{\substack{c \in \mathcal{O}_K \\ (c,1+i)=1}}\mu^2_{[i]}(c) \Phi\left(\frac{N(c)}{X}\right)  \leg{i(1+i)^5c}{\varpi}\sqrt{N(\varpi)}.
\end{split}
\end{align*}

  Our main result evaluates $M_{\Phi}(y,X,\delta)$ in the limit case of $X \rightarrow \infty$.   
\begin{theorem}
\label{quadraticmainthm}
 With the notation as above and the truth of GRH,  we have, for fixed $y\in\mathbb{R}_{>0}$ and $\delta\in(3/4,1)$,  
\begin{align}
\label{quaddensity}
 M_\Phi (y, \delta) \seteq \lim_{X\rightarrow\infty}M_{\Phi}(y,X,\delta) =  \frac{1}{4}\sum_{\substack { l \odd }} \frac {\mu_{[i]}(l)}{N(l^2)} \sum_{\substack{k \in \mathcal{O}_K \\ k \neq 0 }}(-1)^{N(k)}   \widetilde{\Phi}\left(\frac {N(k)}{\sqrt{2y}N(l)}\right),
\end{align}
  where for $t>0$, $\widetilde{\Phi}$ is defined as
\begin{align}
\label{tildeWdef}
   \widetilde{\Phi}(t) =& \int\limits^{\infty}_{-\infty}\int\limits^{\infty}_{-\infty}\Phi(N(x+yi))\widetilde{e}\left(- t(x+yi)\right)\dif x \dif y, \quad \mbox{with} \quad \widetilde{e}(z) =\exp \Big( 2\pi i  \Big( \frac {z}{2i} - \frac {\overline{z}}{2i} \Big) \Big).
\end{align}
  In particular,
\begin{equation}
\label{eq.1leveldensity1}
\lim_{y\rightarrow0^+} M_\Phi(y, \delta)
=0, \quad \text{ and } \quad   \lim_{y\rightarrow \infty} M_{\Phi}(y, \delta) =- \frac{\widetilde{\Phi}(0)}{3\zeta_K(2)}. 
\end{equation}
\end{theorem}

Upon evaluating $\widetilde{\Phi}(t)$ in polar coordinates (see \cite[p. 439]{G&Zhao4}),
\begin{align}
\label{tildeWdefpolr}
   \widetilde{\Phi}(t)= \int\limits^{\infty}_{0}\int\limits^{2\pi}_0\cos (2\pi t r \sin \theta )\Phi(r^2)r  \dif r \dif \theta.
\end{align}

Inserting \eqref{tildeWdefpolr} into \eqref{quaddensity} and a change of variables render that
\begin{align*}
\begin{split}
 M_\Phi (y, \delta) =   \int\limits^{\infty}_{0}\Phi(x)M \Big( \frac yx \Big) \dif x,
\end{split} 
\end{align*}
 where 
\begin{align*}
\begin{split}
 M(x)=\frac{1}{2} \int\limits^{2\pi}_0 \sum_{\substack { l \odd }} \frac {\mu_{[i]}(l)}{N(l^2)} \sum_{\substack{k \in \mathcal{O}_K \\ k \neq 0 }}(-1)^{N(k)} \cos \Big( \frac {\sqrt{2}\pi N(k)}{N(l)\sqrt{x}}\sin \theta \Big) \dif \theta 
\end{split} 
\end{align*}
  is the Zubrilina density for the family $\{ \big ( \frac{i(1+i)^5c}{\cdot} \big): c \in \mathcal O_K, \; (c,1+i) =1, \; c \; \mbox{square-free} \}$ using the terminology introduced in \cite{Sarnak23}. \newline
   
As mentioned earlier, the proof of Theorem \ref{quadraticmainthm} main makes use of the techniques developed in \cite{G&Zhao4} concerning one-level density of the associated family of $L$-functions. We use also the ideas from \cite{LOP25}. As done in \cite[Corollary 1.3]{LOP25}, we note that the expressions in \eqref{tildeWdef} reflect in our case the phase transition range from the corresponding one-level density.  

\section{Preliminaries}
\label{sec 2}

\subsection{Quadratic symbols}
\label{sect: Kronecker}
    For a prime $\varpi \in \mathcal{O}_K$, with $(\varpi, 2)=1$, the quartic symbol is defined for $a \in
\mathcal{O}_K$, $(a, \varpi)=1$ by $\leg{a}{\varpi}_4 \equiv a^{(N(\varpi)-1)/4} \pmod{\varpi}$, with $\leg{a}{\varpi}_4 \in U_K$.  If $\varpi | a$, we define
$\leg{a}{\varpi}_4 =0$.  The quartic character is then extended to any composite $n$ with $(n, 2)=1$ multiplicatively.  Set $\leg{\cdot }{n}_4=1$ for $n \in U_K$. We further define $(\frac{\cdot}{n})=\leg {\cdot}{n}^2_4$  to be the quadratic symbol for these $n$. For $c \in \mathcal O_K$, we write $\chi_c=\leg {c}{\cdot}$ for the quadratic symbol defined above.  \newline

For any $0 \neq q \in \mathcal O_K$, we follow \cite[Section 3.8]{iwakow} to define a Dirichlet character $\chi \pmod {q}$ to be a homomorphism
\begin{align*}
  \chi: \left (\mathcal{O}_K / (q) \right )^{\times}  \rightarrow S^1 :=\{ z \in \mc :  |z|=1 \}.
\end{align*}
We say that $\chi$ is primitive modulo $q$ if it does not factor through $\left (\mathcal{O}_K / (q') \right )^{\times}$ for any divisor $q'$ of $q$ with $N(q')<N(q)$.  Moreover, if $\chi(u)=1$ for any $u \in U_K$, then $\chi$ may be regarded as defined on ideals of $\mathcal O_K$ since every ideal is principal. In this case, we say that this $\chi$ is a Hecke character modulo $q$ of trivial infinite type.  Also, a Hecke character is said to be primitive if it is primitive as a Dirichlet character. 
It is shown in \cite[Section 2.1]{G&Zhao4} that for any $c \in \mathcal O_K$, the symbol $\chi_{i(1+i)^5c}$ is a Hecke character $\pmod {(1+i)^5c}$ of trivial infinite type. Furthermore, if $c$ is square-free, then $\chi_{i(1+i)^5c}$ is non-principal and primitive.

\subsection{Poisson Summation}
The proof of Theorem \ref{quadraticmainthm} requires the following Poisson summation formula from the first display on \cite[p. 442]{G&Zhao4}. 
\begin{lemma}
\label{Poissonsum} With the notation as above, suppose $\Phi$ is a Schwartz class function.  We have, for any primary $l \in \mathcal O_K$ and $X>0$,
\begin{align*}
    \sum_{\substack {a \in \mathcal{O}_K \\ (a, 1+i)=1} } \leg {a}{\varpi}\Phi \left( \frac {N(al^2)}{X} \right)
  & = \frac {X}{2N(l^2)\sqrt{N(\varpi )}} \leg {i(1+i)}{\varpi}\sum_{k \in
   \mathcal{O}_K}(-1)^{N(k)} \leg{k}{\varpi} \widetilde{\Phi}\left(\sqrt{\frac {N(k)X}{2N(l^2\varpi)}}\right), 
\end{align*}
  where $\widetilde{\Phi}$ is defined in \eqref{tildeWdef}.  
\end{lemma}

   We gather \cite[(2.12), (2.14)]{G&Zhao4} to see that for $\Phi$ as given in the statement of Lemma \ref{Poissonsum} that 
for all integers $j \geq 0$, $h \geq 1$ and $t \in \rear$, 
\begin{align}
\label{bounds}
     \widetilde{\Phi}^{(j)}(t) \ll_{j} \frac {1}{(1+|t|)^{h}}. 
\end{align}
   
   We end this section by including a result from \cite[Lemma 3.1]{G&Zhao4}.
 \begin{lemma}
\label{Poissonsumoverk} 
 With the notation as above, we have for any real $A>0$,
\begin{align*}
   \sum_{\substack {k \in \mathcal{O}_K \\ k \neq 0}}(-1)^{N(k)}\widetilde{\Phi}\left(\frac {N(k)}{A}\right)=-\widetilde{\Phi}\left(0\right)+O \left( \frac {1}{A^{1/2}} \right).
\end{align*}
\end{lemma}

\section{Proof of Theorem~\ref{quadraticmainthm}}\label{s:realp}
  
   Let $Z >1$ be a real parameter to be optimized later and decompose $\mu_{[i]}^2(c)$ as
     $\mu_{[i]}^2(c)=M_Z(c)+R_Z(c)$ with
\begin{equation*}
    M_Z(c)=\sum_{\substack {l^2|c \\ N(l) \leq Z \\ l \odd}}\mu_{[i]}(l) \; \quad \mbox{and} \; \quad  R_Z(c)=\sum_{\substack {l^2|c \\ N(l) >
    Z \\ l \odd}}\mu_{[i]}(l).
\end{equation*}
   Since $\Phi$ has compact support, we may define $\beta = \sup_{x \in \mathbb{R}} \{|x| : \Phi(x) > 0\}<\infty$. It follows that, for $Z \in (0, \sqrt{\beta X}]$, 
   \begin{equation} \label{Mdecomp}
   M_{\Phi}(y, X,\delta) = M_{\Phi,Z}(y, X, \delta) + R_{\Phi,Z}(y,X,\delta), 
   \end{equation}
   where
\begin{align}
\label{eq.MAp}
     M_{\Phi,Z}(y, X, \delta) &=\frac{\log X}{X^{1+\delta}} \sum_{\substack{N(\varpi) \in[yX,yX+X^{\delta}] \\ \varpi \odd}} \sum_{\substack{c \in \mathcal{O}_K \\ (c,1+i)=1}}M_Z(c)\Phi\left(\frac{N(c)}{X}\right)  \leg{i(1+i)^5c}{\varpi}\sqrt{N(\varpi)}, 
\end{align}
and $R_{\Phi,Z}(y, X, \delta)$ is the complementary sum defined the same way as $M_{\Phi,Z}(y, X, \delta)$ above, except $R_Z(c)$ replaces $M_Z(c)$.

\subsection{Bounding $R_{\Phi,Z}(y,X, \delta)$}
  First we observe that for a given $c \in \mathcal{O}_K$, for any $\varepsilon>0$, 
\begin{equation}\label{eq.sm}
    \sum_{\substack{l^2 \mid c , \ N(l) > Z \\ l \odd}} \mu_{[i]}(l)  \ll \sum_{l \mid c} 1 \ll N(c)^\varepsilon,
\end{equation}
  where the last estimation above follows by arguing similar to the proof of \cite[Theorem 2.11]{MVa1}. \newline
  
Now writing $c = l^2b$ with $b \odd$ and applying \eqref{eq.sm} lead to
\begin{equation*}
 R_{\Phi,Z}(y, X, \delta) \ll \frac{\log X}{X^{1/2+\delta-\varepsilon}}  \sum_{ N(l) \in (Z,\sqrt{\beta X}]} \sum_{N(b) \leq \frac{\beta X}{N(l)^2} } \Phi\left(\frac{N(l^2b)}{X}\right)\Bigg|\sum_{\substack{N(\varpi)\in[yX,yX+X^{\delta}] \\ \varpi \odd }}\leg{i(1+i)^5l^2b}{\varpi} \sqrt{\frac {N(\varpi)}{X}} \Bigg|. 
\end{equation*}
It follows from \cite[Lemma 2.5]{G&Zhao4} that under GRH, for any non-principal Hecke character $\chi$ of trivial infinite type with modulus dividing $n$
\begin{align} 
\label{lem2.7eq}
 \sum_{\substack {N(\varpi) \leq x \\ \varpi \odd}} \chi (\varpi) \log N(\varpi)
\ll \sqrt{x} \log^{3} x \log (N(n)+2), \quad \mbox{for} \; x \geq 1.
\end{align}  

Our discussions in Section \ref{sect: Kronecker} gives that for any odd $c \in \mathcal{O}_K$, the symbol $\leg {i(1+i)^5c}{\cdot}$ is a Hecke character of trivial infinite type whose conductor dividing $c$.  Thus, from \eqref{lem2.7eq} and partial summation that, for odd $c \in \mathcal{O}_K$ and $x \geq 2$, under GRH, 
\begin{align}
\label{sumchiest}
\begin{split}    
  S(x, \chi_{i(1+i)^5c}) :=\sum_{N(\varpi) \leq x} \chi_{i(1+i)^5c}(\varpi) \ll x^{1/2}\log (N(c)x). 
\end{split}
\end{align}
The above, with partial summation, yeilds that, under GRH, for any fixed $y\in\mathbb{R}_{>0}$, $\varepsilon>0$ and sufficiently large $X$, 
\begin{equation*}
\sum_{\substack{\varpi \in[yX,yX+X^{\delta}] \\ \varpi \odd }} \leg{i(1+i)^5l^2b}{\varpi} \frac {\sqrt{N(\varpi)}}{\sqrt{X}} \ll y^{1+\varepsilon}X^{1/2+\varepsilon}.
\end{equation*}
   
Now setting $Z = X^{1/4}$ and proceeding as in Section 4.1 of \cite{LOP25}, we get that for $\varepsilon_0=(\delta-3/4)/5$,
\begin{equation}\label{eq.RllXe}
    R_{\Phi,Z}(y,X,\delta) \ll y^{1+\varepsilon_0}X^{-\varepsilon_0}.
\end{equation}
Using \eqref{eq.RllXe} and \eqref{Mdecomp}, we obtain that
\begin{equation*}
\lim_{X\rightarrow\infty}M_{\Phi}(y,X,\delta)=\lim_{X\rightarrow\infty}M_{\Phi,Z}(y,X,\delta).
\end{equation*}

\subsection{Evaluation of $M_{\Phi,Z}(y,X, \delta)$}
\label{eq.anoMPA}

From \eqref{eq.MAp},  we write
\begin{equation} \label{3.0}
 M_{\Phi,Z}(y,X, \delta) = \frac{\log X}{X^{1+\delta}} \sum_{\substack{N(\varpi) \in[yX,yX+X^{\delta}] \\ \varpi \odd}} \sqrt{N(\varpi)} \leg{i(1+i)}{\varpi} \sum_{\substack {N(l) \leq Z \\ l \odd}} \mu_{[i]}(l)\leg {l^2}{\varpi}  \sum_{\substack {a \in \mathcal{O}_K \\ (a, 1+i)=1} } \leg {a}{\varpi}\Phi \left( \frac {N(al^2)}{X} \right).
\end{equation}
 
Now Lemma~\ref{Poissonsum}, applied to the inner-most sum above, recasts $M_{\Phi,Z}(y,X, \delta)$ as
\begin{align*}
\begin{split}
   M_{\Phi,Z}(y,X, \delta)   =& \frac{\log X}{2X^\delta} \sum_{\substack{N(\varpi) \in[yX,yX+X^{\delta}] \\ \varpi \odd}} \sum_{\substack {N(l) \leq Z \\ l \odd }} \frac {\mu_{[i]}(l)}{N(l^2)} \sum_{\substack{k \in
   \mathcal{O}_K \\ k \neq 0}}(-1)^{N(k)}  \leg {k}{\varpi}\widetilde{\Phi}\left(\sqrt{\frac {N(k)X}{2N(l^2\varpi)}}\right),
   \end{split}
\end{align*}
observing that the contribution from the terms with $k=0$ is zero.  Moreover, due to the presence of $\big( \frac{l^2}{\varpi} \big)$ in \eqref{3.0}, we also have the condition $(l, \varpi) = 1$ in the above sum.  But as $N(l) \leq Z = X^{1/4} \ll N(\varpi)$, that co-primality is automatically satisfied for sufficiently large $X$. \newline

We write $\square$ for a perfect square in $\mathcal O_K$.  The part of the sum $M_{\Phi,Z}(y,X, \delta)$
corresponding to the contribution of $k \neq 0, \square$ can be written as $(\log X)R/(2X^{\delta})$, where
\begin{align*}
R:= \sum_{\substack {N(l) \leq Z \\ l \odd }} \frac {\mu_{[i]}(l)}{N(l^2)} \sum_{\substack {k \in
   \mathcal{O}_K \\ k \neq 0, \square}} & (-1)^{N(k)}  \sum_{\substack{N(\varpi) \in[yX,yX+X^{\delta}] \\ \varpi \odd}}  \leg {k}{\varpi} \widetilde{\Phi}\left(\sqrt{\frac {N(k)X}{2N(l^2\varpi)}}\right).
\end{align*}

   Recall from our discussions in Section \ref{sect: Kronecker} that we denote $\chi_{k}$ for the quadratic symbol $\leg {k}{\cdot}$ and when $k$ is not a square, $\chi_{k}$ can be regarded as a non-principle Hecke character modulo $k$ of trivial infinite type.  Now \eqref{sumchiest} and partial summation reveal that
\begin{align*}
\begin{split}
 \sum_{\substack{N(\varpi) \in[yX,yX+X^{\delta}] \\ \varpi \odd}} & \leg {k}{\varpi} \widetilde{\Phi}\left(\sqrt{\frac {N(k)X}{2N(l^2\varpi)}}\right) =\int\limits^{yX+X^{\delta}}_{yX}\widetilde{\Phi}\left(\sqrt{\frac {N(k)X}{2N(l^2)V}}\right) \dif S \left( V, \chi_{k}\right)  \\
& \hspace*{1cm} \ll (yX)^{1/2}\log \left( yX(N(k)+2) \right) \left( 
   \left | \widetilde{\Phi}\left(\sqrt{\frac {N(k)}{2N(l^2)y}}\right) \right |  +   \left |\widetilde{\Phi}\left(\sqrt{\frac {N(k)X}{2N(l^2)(yX+X^{\delta})}}\right) \right | \right.\\
  & \hspace*{7cm} \left.+ \int\limits^{yX+X^{\delta}}_{yX}\sqrt{\frac {N(k)X}{N(l^2)}}\frac
   {1}{V^{3/2}} \left |\widetilde{\Phi}'\left(\sqrt{\frac {N(k)X}{2N(l^2)V}}\right) \right | \dif V \right).
\end{split}
\end{align*}

This gives raise to
\begin{align*}
  R  \ll \sum_{\substack {N(l) \leq Z \\ l \odd}}\frac {1}{N(l^2)} \left( R_1+R_2+R_3 \right), \quad \mbox{where} \quad   R_1  =  (yX)^{1/2}\sum_{\substack {k \in
   \mathcal{O}_K \\ k \neq 0}}\log \left( yX(N(k)+2) \right) \left |\widetilde{\Phi}\left(\sqrt{\frac {N(k)}{2N(l^2)y}}\right) \right | ,
\end{align*}
\[ R_2  = (yX)^{1/2}\sum_{\substack {k \in  \mathcal{O}_K \\ k \neq 0}}\log \left( yX(N(k)+2) \right) \left |\widetilde{\Phi}\left(\sqrt{\frac {N(k)X}{2N(l^2)(yX+X^{\delta})}}\right) \right |, \]
and
\begin{equation*}
   R_3 =  y^{1/2}X\int\limits^{yX+X^{\delta}}_{yX}\frac
   {1}{N(l)V^{3/2}}  \sum_{\substack {k \in
   \mathcal{O}_K \\ k \neq 0}}\log(yX(N(k)+2)) \sqrt{N(k)} \left |\widetilde{\Phi}'\left(\sqrt{\frac {N(k)X}{2N(l^2)V}}\right) \right |\dif V.
\end{equation*}

Now \eqref{bounds} with $h=3$ implies that
\begin{align*}
  R_1 \ll & (yX)^{1/2}\sum_{\substack {k \in
   \mathcal{O}_K \\ 1 \leq N(k) \leq yN(l^2)}}\log \left( yX(N(k)+2) \right)+(yX)^{1/2}\sum_{\substack {k \in
   \mathcal{O}_K \\  N(k) > yN(l^2)}}\log \left( yX(N(k)+2) \right)\left ( \frac {N(l^2)y}{N(k)}\right )^{3/2} \\
 \ll&  (yX)^{1/2+\varepsilon}yN(l^2).
\end{align*}
  Similarly,
\begin{align*}
  R_2  \ll&  (yX)^{1/2+\varepsilon}\frac {yX+X^{\delta}}{X}N(l^2) \ll (yX)^{1/2+\varepsilon}yN(l^2), \quad \mbox{and} \\
  R_3  \ll &  \frac {X^{1/2}(yX)^{1/2+\varepsilon}}{N(l)} \int\limits^{yX+X^{\delta}}_{yX}\frac
   {1}{V^{3/2}}  \left(\frac {N(l^2)V}{X}\right)^{3/2} \dif V \ll \frac {X^{1/2}(yX)^{1/2+\varepsilon}N(l)^2}{X^{1-\delta}}. 
\end{align*}

   These estimates give
\begin{equation*}
  R \ll  (yX)^{1/2+\varepsilon}yZ.
\end{equation*}
and thus we conclude that the contribution of $k \neq 0$, $\square$ to $M_{\Phi,Z}(y,X, \delta)$ is, for any $\varepsilon>0$, 
\begin{equation*}
 \ll  \frac {(yX)^{1/2+\varepsilon}yZ}{X^{\delta}}.
\end{equation*}

Recall from the discussion before \eqref{eq.RllXe} that $Z = X^{1/4}$ and $\delta>3/4$.  Hence the above vanishes in the limit as $X\rightarrow\infty$ and it suffices to evaluate
\begin{equation*}
  \lim_{X\rightarrow\infty}\frac{\log X}{2X^\delta}\sum_{\substack{N(\varpi) \in[yX,yX+X^{\delta}] \\ \varpi \odd}}\sum_{\substack {N(l) \leq Z \\ l \odd }} \frac {\mu_{[i]}(l)}{N(l^2)} \sum_{\substack{k \in
   \mathcal{O}_K \\ k \neq 0 \\ k =\square }}(-1)^{N(k)}    \leg {k}{\varpi}\widetilde{\Phi}\left(\sqrt{\frac {N(k)X}{2N(l^2\varpi)}}\right).
\end{equation*}

Changing variables $k \mapsto k^2$ in the inner-most sum above and noting that $k^2=k_1^2$ if and only if $k =\pm k_1$, we deduce
\begin{align}
\label{MS}
\begin{split}
  \lim_{X\rightarrow\infty}M_{\Phi,Z}(y,X,\delta) = &  \lim_{X\rightarrow\infty}\frac{\log X}{4X^\delta}\sum_{\substack{N(\varpi) \in[yX,yX+X^{\delta}] \\ \varpi \odd}}\sum_{\substack {N(l) \leq Z \\ l \odd }} \frac {\mu_{[i]}(l)}{N(l^2)} \sum_{\substack{k \in
   \mathcal{O}_K \\ k \neq 0 \\ (k, \varpi)=1 }}(-1)^{N(k)}   \widetilde{\Phi}\left(N(k)\sqrt{\frac {X}{2N(l^2\varpi)}}\right) \\
 =: &\lim_{X\rightarrow\infty}S-\lim_{X\rightarrow\infty}S', 
\end{split}
\end{align}
   where $S$ discards the condition $(k, \varpi)=1$ and $S'$ has its inner-most sum over multiples of $\varpi$, i.e. 
\begin{align*}
\begin{split}
  S =& \frac{\log X}{4X^\delta}\sum_{\substack{N(\varpi) \in[yX,yX+X^{\delta}] \\ \varpi \odd}}\sum_{\substack {N(l) \leq Z \\ l \odd }} \frac {\mu_{[i]}(l)}{N(l^2)} \sum_{\substack{k \in
   \mathcal{O}_K \\ k \neq 0  }}(-1)^{N(k)}   \widetilde{\Phi}\left(N(k)\sqrt{\frac {X}{2N(l^2\varpi)}}\right), \quad \mbox{and} \\
   S' =&\frac{\log X}{4X^\delta}\sum_{\substack{N(\varpi) \in[yX,yX+X^{\delta}] \\ \varpi \odd}}\sum_{\substack {N(l) \leq Z \\ l \odd }} \frac {\mu_{[i]}(l)}{N(l^2)} \sum_{\substack{k \in
   \mathcal{O}_K \\ k \neq 0  }}(-1)^{N(k)}   \widetilde{\Phi}\left(N(k)\sqrt{\frac {XN(\varpi)}{2N(l^2)}}\right).
\end{split}
\end{align*}

    We now estimate $S'$ by applying \eqref{bounds} with $h=2$, getting
\begin{align*}
    \sum_{\substack {k \in
  \mathcal{O}_K \\ k \neq 0}}(-1)^{N(k)}\widetilde{\Phi}\left(N(k) \sqrt{\frac {X N(\varpi) }{2N(l^2)}}\right) \ll
   \sum_{\substack {k \in
   \mathcal{O}_K  \\ N(k) \geq 1}} \frac {N(l^2)}{N^2(k) X N(\varpi)} \ll \frac {N(l^2)}{X N(\varpi)}.
\end{align*}

From this, emerges the bound
\begin{align*}
   S' \ll \frac{(\log X)^2Z}{X^{1+\delta}}.
\end{align*}
  As $Z = X^{1/4}$, the above, applied to \eqref{MS}, leads to  
\begin{align*}
\begin{split}
  \lim_{X\rightarrow\infty}M_{\Phi,Z}(y,X,\delta) = & \lim_{X\rightarrow\infty}S.
\end{split}
\end{align*}

To evaluate $S$, we rewrite the inner-most sum of $S$ as
\[ -\widetilde{\Phi}(0)+T(X, N(l), N(\varpi)), \quad \mbox{where} \quad
T(X, N(l), N(\varpi)) =  \sum_{\substack{k \in \mathcal{O}_K  }}(-1)^{N(k)}   \widetilde{\Phi}\left(N(k)\sqrt{\frac {X}{2N(l^2\varpi)}}\right) . \]
By Lemma \ref{Poissonsumoverk},
\begin{align}
\label{Tbound}
\begin{split}
  T(X, N(l), N(\varpi)) \ll \Big (\frac {X}{N(l^2\varpi)} \Big )^{1/4}. 
\end{split}
\end{align} 

Now we get
\begin{align}
\label{Sexp}
\begin{split}
  S = & \frac{\log X}{4X^\delta}\sum_{\substack{N(\varpi) \in[yX,yX+X^{\delta}] \\ \varpi \odd}} \sum_{\substack {l \odd }} \frac {\mu_{[i]}(l)}{N(l^2)} (-\widetilde{\Phi}(0)+T(X, N(l), N(\varpi))) \\
& \hspace*{2cm} - \frac{\log X}{4X^\delta}\sum_{\substack{N(\varpi) \in[yX,yX+X^{\delta}] \\ \varpi \odd}}\sum_{\substack {N(l) >Z  \\ l \odd }} \frac {\mu_{[i]}(l)}{N(l^2)} (-\widetilde{\Phi}(0)+T(X, N(l), N(\varpi))) \\
=&  \frac{\log X}{4X^\delta}\sum_{\substack{N(\varpi) \in[yX,yX+X^{\delta}] \\ \varpi \odd}}\sum_{\substack {l \odd }} \frac {\mu_{[i]}(l)}{N(l^2)} (-\widetilde{\Phi}(0)+T(X, N(l), N(\varpi))) \\
& \hspace*{2cm}+ \frac{\log X}{4X^\delta}\sum_{\substack{N(\varpi) \in[yX,yX+X^{\delta}] \\ \varpi \odd}} O \Big(\frac {1}{Z}+\frac {X^{1/4}}{Z^{3/2}N(\varpi)^{1/4}} \Big).
\end{split}
\end{align}

  We note that the prime ideal theorem under GRH asserts that (see \cite{LO77}) for $x >2$, 
\begin{align*}
\begin{split}
  \sum_{\substack{\varpi  \\ N(\varpi)  \leq x  }}\log N(\varpi) =x+O(x^{1/2}\log^2 x). 
\end{split}
\end{align*}
  We apply the above and partial summation to see that for any $y\in\mathbb{R}_{>0}$ and $\delta\in(3/4,1)$, we have
\begin{align}
\label{sumvarpi}
\begin{split}
  \sum_{\substack{N(\varpi) \in[yX,yX+X^{\delta}] \\ \varpi \odd}}1=& \frac {X^{\delta}}{\log X}+o\Big( \frac {X^{\delta}}{\log X} \Big), \quad \mbox{and} \\
  \sum_{\substack{N(\varpi) \in[yX,yX+X^{\delta}] \\ \varpi \odd}}N(\varpi)^{-1/4} \ll & (yX)^{-1/4}\sum_{\substack{N(\varpi) \in[yX,yX+X^{\delta}] \\ \varpi \odd}}1 \ll (yX)^{-1/4}\frac {X^{\delta}}{\log X}. 
\end{split}
\end{align}

   Recall that $Z=X^{1/4}$. Then the bounds in \eqref{sumvarpi} imply that
\begin{align*}
\begin{split}
 \lim_{X \rightarrow \infty} \frac{\log X}{4X^\delta}\sum_{\substack{N(\varpi) \in[yX,yX+X^{\delta}] \\ \varpi \odd}} \Big(\frac {1}{Z}+\frac {X^{1/4}}{Z^{3/2}N(\varpi)^{1/4}}\Big)=0.
\end{split}
\end{align*}
 
   We arrive at, using the above and \eqref{Sexp}, 
\begin{align*}
\begin{split}
 \lim_{X \rightarrow \infty} S 
=&  \lim_{X \rightarrow \infty} \frac{\log X}{4X^\delta}\sum_{\substack{N(\varpi) \in[yX,yX+X^{\delta}] \\ \varpi \odd}}\sum_{\substack {l \odd }} \frac {\mu_{[i]}(l)}{N(l^2)} (-\widetilde{\Phi}(0)+T(X, N(l), N(\varpi))). 
\end{split}
\end{align*}

  We apply \eqref{Tbound} and the Cauchy criterion for uniform convergence to see that the series
\begin{align*}
\begin{split}
 \sum_{\substack {l \odd }} \frac {\mu_{[i]}(l)}{N(l^2)} (-\widetilde{\Phi}(0)+T(X, N(l), N(\varpi)))
\end{split}
\end{align*} 
  converges uniformly and hence represents a continuous function of $N(\varpi)$. We then argue in a way similar to the proof of \cite[Lemma 2.9]{LOP25}, using \eqref{sumvarpi}.  This leads to, under GRH,
\begin{align}
\label{Sexp1}
\begin{split}
 M_\Phi (y, \delta) = \lim_{X\rightarrow\infty}M_{\Phi}(y,X,\delta) = \lim_{X \rightarrow \infty} S 
=&  \frac{1}{4}\sum_{\substack {l \odd }} \frac {\mu_{[i]}(l)}{N(l^2)}\sum_{\substack{k \in
   \mathcal{O}_K \\ k \neq 0  }}(-1)^{N(k)}   \widetilde{\Phi}\left(N(k)\sqrt{\frac {1}{2yN(l^2)}}\right)
\end{split}
\end{align}
and establishes \eqref{quaddensity}.   

\subsection{Evaluation of $\lim_{y \rightarrow \infty}M_\Phi (y, \delta)$}
 
  To evaluate the second limit in \eqref{eq.1leveldensity1}, we observe that similar to our discussions above,  the function $M_\Phi (y, \delta)$ is a continuous function of $y$. Moreover,
\begin{align}
\label{Sexp2}
\begin{split}
 \sum_{\substack{k \in
   \mathcal{O}_K \\ k \neq 0  }}(-1)^{N(k)}   \widetilde{\Phi}\left(N(k)\sqrt{\frac {X}{2yN(l^2)}}\right)=-\widetilde{\Phi}(0)+T(1, N(l), y). 
\end{split}
\end{align}
By \eqref{Tbound},
\begin{align*}
\begin{split}
  T(1, N(l), y) \ll \Big (\frac {1}{N(l^2)y} \Big )^{1/4}. 
\end{split}
\end{align*} 

It follows that
\begin{align*}
\begin{split}
 \lim_{y \rightarrow \infty}M_\Phi (y, \delta) = - \frac{1}{4}\widetilde{\Phi}(0)\sum_{\substack {l \odd }} \frac {\mu_{[i]}(l)}{N(l^2)}= - \frac{\widetilde{\Phi}(0)}{3\zeta_K(2)}. 
\end{split}
\end{align*}
  This establishes the second limit given in \eqref{eq.1leveldensity1}.

\subsection{Evaluation of $\lim_{y \rightarrow 0^+}M_\Phi (y, \delta)$}

  To compute $\lim_{y \rightarrow 0^+}M_\Phi (y, \delta)$, \cite[(3.8)]{G&Zhao4} gives that
\begin{align}
\label{Sexp4}
\begin{split}
 T(1, N(l), y)=\sqrt{2yN(l^2)}\sum_{\substack{ j \in
   \mathcal{O}_K \\ (j, 1+i)=1}}\breve{\Phi}\left(\sqrt{\frac {N(j)\sqrt{2yN(l^2)}}{2}}\right),
\end{split}
\end{align}
 where 
\begin{align*}
   \breve{\Phi}(t) =\int\limits^{\infty}_{-\infty}\int\limits^{\infty}_{-\infty}\widetilde{\Phi}(N(u+vi))\widetilde{e}\left(- t(u+vi)\right)\dif u \dif v, \quad t \geq 0.
\end{align*}

Mellin inversion leads to
\begin{align}
\label{Phix}
 \begin{split}
  \breve{\Phi}\left(x \right)=\int\limits_{(c)}\hat{\Phi}(s)x^{-s} \dif s,
 \end{split}
 \end{align}
   where we choose $c$ with $0<c<1$ and
\begin{align*}
 \begin{split}
  \hat{\Phi}(s)=\int\limits^{\infty}_0 \breve{\Phi}(x)x^s\frac {\dif x}{x}. 
 \end{split}
 \end{align*}
Similar to \eqref{tildeWdefpolr}, 
\begin{equation*}
\begin{split}
    \breve{\Phi}(x) =\int\limits_{\mr^2}\cos (2\pi x v)\widetilde{\Phi}(u^2+v^2) \dif u \dif v & =4\int\limits^{\infty}_{0}\int\limits^{\pi/2}_0\cos (2\pi x r \sin \theta )\widetilde{\Phi}(r^2)r \dif r \dif \theta.
    \end{split}
\end{equation*}
Hence, for $0<\Re(s)<1$, 
\begin{align*}
 \begin{split}
  \hat{\Phi}(s)= \int^{\infty}_0 \breve{\Phi}(x)x^s\frac {dx}{x}= & 4\int\limits^{\infty}_{0}\int\limits^{\pi/2}_0\int\limits^{\infty}_0 \cos (2\pi x r \sin \theta )x^{s}\frac {dx}{x}\widetilde{\Phi}(r^2)r \dif r \dif \theta \\
=& \frac{4}{(2\pi)^s}\int\limits^{\infty}_{0}\int\limits^{\pi/2}_0\int\limits^{\infty}_0 \cos (x)x^{s}\frac {\dif x}{x}\widetilde{\Phi}(r^2)r^{1-s}(\sin \theta)^{-s} \dif r \dif \theta.
 \end{split}
 \end{align*}
   Using the formula (see \cite[17.43.3]{GR})
\begin{align*}
 \begin{split}
  \int\limits^{\infty}_0 \cos (x)x^{s}\frac {\dif x}{x}=\Gamma (s)\cos \Big( \frac {\pi s}{2} \Big) \quad \mbox{for} \quad 0 < \Re (s) < 1,
 \end{split}
 \end{align*}
  we get 
\begin{align}
\label{hatphi0}
 \begin{split}
  \hat{\Phi}(s)= \frac{4\Gamma (s)}{(2 \pi)^s} \cos \Big(\frac {\pi s}{2} \Big)\int\limits^{\infty}_{0}\int\limits^{\pi/2}_0\widetilde{\Phi}(r^2)r^{1-s}(\sin \theta)^{-s} \dif r \dif \theta.
 \end{split}
 \end{align}

Moreover, by \cite[3.621.5]{GR} and \cite[8.384.1]{GR}, we have, for $0<\Re(s)<1$,  
\begin{align}
\label{hatphi1}
 \begin{split}
  \int\limits^{\pi/2}_0(\sin \theta)^{-s}\dif \theta= \frac {\Gamma(\tfrac12)\Gamma(\tfrac{1-s}{2})}{2\Gamma(\tfrac{2-s}{2})}.
 \end{split}
 \end{align}

   Also, integration by parts renders that, for $0<\Re(s)<1$,  
\begin{align}
\label{hatphi2}
 \begin{split}
  \int\limits^{\infty}_{0}\widetilde{\Phi}(r^2)r^{1-s}\dif r =\frac 12\int\limits^{\infty}_{0}\widetilde{\Phi}(r)r^{-s/2}\dif r=-\frac 1{2-s}\int\limits^{\infty}_{0}\widetilde{\Phi}'(r)r^{(2-s)/2}\dif r.
 \end{split}
 \end{align}

From \eqref{hatphi0}--\eqref{hatphi2} that 
\begin{align}
\label{hatphi3}
 \begin{split}
  \hat{\Phi}(s)=& -\frac {\Gamma (s)\Gamma(\frac 12)\cos (\frac {\pi s}{2})\Gamma(\frac {1-s}2)}{(2\pi)^s(\frac {2-s}{2})\Gamma(\frac {2-s}2)}\int\limits^{\infty}_{0}\widetilde{\Phi}'(r)r^{(2-s)/2}\dif r.
 \end{split}
 \end{align}

   Note that the above expression is initially valid for  $0<\Re(s)<1$ and remains holomorphic for $0<\Re(s)<5/2$, since the pole of $\Gamma(\frac {1-s}2)$ at $s=1$ (resp. $(2-s)^{-1}$ at $s=2$) is cancelled by the zero of $\cos (\frac {\pi s}{2})$ (resp. $\Gamma^{-1}(\frac {2-s}2)$). Furthermore, by Stirling's formula (see \cite[(5.113)]{iwakow}), for $0<\Re(s)<5/2$,
\begin{align}
\label{Stirling}
 \begin{split}
 \Gamma(s)\cos \Big(\frac {\pi s}{2} \Big) \ll (1+|s|)^{\Re(s)-1/2} \quad \mbox{and} \quad \frac {\Gamma(\frac {1-s}2)}{(\frac {2-s}{2})\Gamma(\frac {2-s}2)} \ll (1+|s|)^{-3/2}.
 \end{split}
 \end{align}
    Moreover, integration by parts and \eqref{bounds} give that, for any integer $E \geq 1$ and $0<\Re(s)<5/2$,
\begin{align}
\label{intbound}
 \begin{split}
  \int\limits^{\infty}_{0}\widetilde{\Phi}'(r)r^{(2-s)/2}\dif r \ll (1+|s|)^{-E}.
 \end{split}
 \end{align}

   It follows from \eqref{Stirling} and \eqref{intbound} that, for  any integer $E \geq 1$ and $0<\Re(s)<5/2$,
\begin{align}
\label{hatphibound}
 \begin{split}
  \hat{\Phi}(s)\ll (1+|s|)^{-E}.
 \end{split}
 \end{align}

We apply \eqref{Stirling} and \eqref{intbound} to see that shifting the line of integration in \eqref{Phix} to $c=5/2$ encounters no pole and the resulting integral remains convergent. Thus
\begin{align*}
 \begin{split}
  \breve{\Phi}\left(x \right)=\int\limits_{(5/2)}\hat{\Phi}(s)x^{-s} \dif s,
 \end{split}
 \end{align*} 
   
The above, together with \eqref{Sexp1}--\eqref{Sexp4}, allows us to deduce
\begin{align}
\label{eq:m_phi_divsum}
\begin{split}
    M_\Phi(y, \delta) =&  - \frac{\widetilde{\Phi}(0)}{3\zeta_K(2)} + \frac{\sqrt{2y}}{4}\sum_{\substack {l \odd }} \frac {\mu_{[i]}(l)}{N(l)}\sum_{\substack{ j \in
   \mathcal{O}_K \\ (j, 1+i)=1}}\breve{\Phi}\left(\sqrt{N(j)N(l)\sqrt{\frac {y}{2}}}\right)\\
    =& - \frac{\widetilde{\Phi}(0)}{3\zeta_K(2)} + \frac{\sqrt{2y}}{4}\sum_{\substack{ j \in
   \mathcal{O}_K \\ (j, 1+i)=1}}\Big (\sum_{\substack {l \mid j \\ l \odd }} \frac {\mu_{[i]}(l)}{N(l)}\Big )\int\limits_{(5/2)}\hat{\Phi}(s) \Big(N(j)\sqrt{\frac {y}{2}}\Big)^{-s/2} \dif s \\
 =& - \frac{\widetilde{\Phi}(0)}{3\zeta_K(2)} + \sqrt{2y}\int\limits_{(5/2)}\hat{\Phi}(s)B\Big(\frac s2\Big) \Big(\sqrt{\frac {y}{2}} \Big)^{-s/2} \dif s, 
 \end{split}   
 \end{align}
  where we define, for $\Re(s)>1$, 
\begin{align*}
 \begin{split}
      B(s) &=\sum_{n \odd} \Big (\sum_{\substack {l \mid n \\ l \odd }} \frac {\mu_{[i]}(l)}{N(l)}\Big ) N(n)^{-s} .
 \end{split}
 \end{align*}   
   Similar to \cite[(5.4)]{LOP25}, for $\mathrm{Re}(s)>1$, we have
 \begin{align}
 \label{eq.Bszeta}
\begin{split}
      B(s) &=\frac {1-2^{-s}}{1-2^{-s-1}} \cdot \frac{\zeta_K(s)}{\zeta_K(s+1)} .
 \end{split}
 \end{align}
  The above implies that $B(s)$  has meromorphic continuation to $\mathrm{Re}(s) > 0$ with a simple pole at $s=1$ and (recall that the residue of
$\zeta_K(s)$ at $s=1$ equals $\pi/4$)
 \begin{equation}
\label{eq.Rs1B}
     \mathrm{Res}_{s=1}B(s) = \frac {\pi}{6}  \cdot \frac 1 {\zeta_K(2)}.
 \end{equation}

    Also, it is shown in \cite[p. 17]{iwakow} that
\begin{align}
\label{zetaKdecomp}
\zeta_K(s)=\zeta(s)L(s,\psi_4),
\end{align}
   where $\psi_4$ the only primitive Dirichlet character modulo $4$, i.e. $\psi_4(n) =(-1)^{(n-1)/2}$ if $(n,2)=1$, $\psi_4(n)=0$ if $2|n$. \newline
  
   We then deduce from \cite[(5.20)]{iwakow} that for $\Re(s) \geq 1/2, |s-1|>1/2$ and any $\varepsilon>0$, 
\begin{align}
\label{zetaKbound}
\zeta_K(s) \ll (1+|s|)^{1-\Re(s)+\varepsilon},
\end{align}

A Vinogradov-Korobov type zero-free region t is established in \cite{Sokolovskii68} for $\zeta_K(s)$.  Namely, for $s=\sigma+it$ with $\sigma$, $t \in \mr$, there are suitable positive constants $c$ and $t_0$, such that $\zeta_K(s) \neq 0$ in the region
\begin{align*}
 \sigma \geq 1-c(\log |t|)^{-2/3}(\log \log |t|)^{-1/3}, \quad |t| \geq t_0. 
\end{align*}  
   In view of \eqref{zetaKdecomp}, the above now implies that neither $\zeta(s)$ nor $L(s,\psi_4)$ vanishes in the above region.  From this and \cite[Theorems 6.7, 11.4]{MVa1}, we have for $\Re(s) \geq 1$, 
\begin{align}
\label{zetaKinvbound}
\zeta^{-1}_K(s) \ll \log^2 (2+|s|). 
\end{align}
  
Now, using \eqref{eq.Bszeta}, \eqref{zetaKbound} and \eqref{zetaKinvbound}, we get for $\Re(s) \geq 1/2$, 
\begin{align}
 \label{eq.Bsbound}
\begin{split}
      B(s) \ll & (1+|s|)^2.
 \end{split}
 \end{align}  
    We now shift the line of integration of the last expression in \eqref{eq:m_phi_divsum} to $\Re(s)=3/2$ to see from \eqref{eq.Rs1B}, \eqref{hatphibound}, \eqref{eq.Bsbound} that
\begin{align}
\label{eq:m_phi_divsum1}
\begin{split}
    M_\Phi(y, \delta) 
 =& - \frac{\widetilde{\Phi}(0)}{3\zeta_K(2)} + \sqrt{2y}\text{Res}_{s=2}\Big (\hat{\Phi}(s)B\Big(\frac s2\Big)(\sqrt{\frac {y}{2}})^{-s/2}\Big )+O(y^{1/8}) \\
=& - \frac{\widetilde{\Phi}(0)}{3\zeta_K(2)} + \sqrt{2y}\frac {2\pi}{6}  \cdot \frac{\pi}{3\zeta_K(2)}\hat{\Phi}(2)\Big(\sqrt{\frac {y}{2}}\Big)^{-1}+O(y^{1/8}) = - \frac{\widetilde{\Phi}(0)}{3\zeta_K(2)} + \frac {2\pi\hat{\Phi}(2)}{3}  \cdot \frac 1 {\zeta_K(2)}+O(y^{1/8}).  
 \end{split}   
 \end{align}
 
   It remains to evaluate $\hat{\Phi}(2)$. For this, note that (see \cite[(C.8)]{MVa1}) $\Gamma(s)$ has a simple pole at $s=-n$ for every $n \in \{0,1,2, \cdots \}$ with residue $(-1)^n/n!$. It follows from this and \eqref{hatphi3} that
 \begin{align}
\label{hatphi3-2}
 \begin{split}
  \hat{\Phi}(2)=& -(2\pi)^{-2}\Gamma (2)\Gamma(\tfrac 12)\cos (\pi)\Gamma(-\tfrac 12)\int\limits^{\infty}_{0}\widetilde{\Phi}'(r)\dif r=-(2\pi)^{-2}\Gamma(\tfrac 12)\Gamma(-\tfrac 12)\widetilde{\Phi}(0).
 \end{split}
 \end{align}

 The well-known relation $s\Gamma(s)=\Gamma(s+1)$ (see \cite[Chap. 10 (3)]{Da}) implies that $\Gamma(-\tfrac12)=-2\Gamma(\tfrac12)$.   Moreover, (see \cite[(C.7)]{MVa1}) $\Gamma(\tfrac12)=\sqrt{\pi}$. It follows from this and \eqref{hatphi3-2} that
 \begin{align*}
 \begin{split}
  \hat{\Phi}(2)=& (2\pi)^{-1}\widetilde{\Phi}(0).
 \end{split}
 \end{align*}

Inserting the above into \eqref{eq:m_phi_divsum1} yields 
 \begin{align*}
\begin{split}
    M_\Phi(y, \delta) =O(y^{1/8}),
 \end{split}   
 \end{align*}
from which we conclude
 \begin{equation*}
     \lim_{y \to 0^+} M_\Phi(y,\delta) = 0.
 \end{equation*}
   This establishes the first limit given in \eqref{eq.1leveldensity1} and thus completes the proof of Theorem \ref{quadraticmainthm}. \newline
  
\noindent{\bf Acknowledgments.} P. G. is supported in part by NSFC grant 12471003 and L. Z. by the FRG Grant PS71536 at the University of New South Wales.

\bibliography{biblio}
\bibliographystyle{amsxport}

\end{document}